\documentclass[11pt]{article}

\usepackage{amsfonts}
\usepackage{amssymb,amsmath} 
\usepackage{psfrag,graphicx}

\newcommand{\nicefrac}[2]
{\leavevmode \kern.1em\raise.5ex\hbox{\the\scriptfont0 #1}
             \kern-.1em/\kern-.15em\lower.25ex
             \hbox{\the\scriptfont0 #2}}

\newtheorem{Theo}{Theorem}{\alph{enumi}}
\newenvironment{theorem}{\begin{Theo}\hspace{-0.2cm}: }{\end{Theo}}

\newtheorem{Pro}{Proposition}{\alph{enumi}}
\newenvironment{proposition}{\begin{Pro}\hspace{-0.2cm}: }{\end{Pro}}

\newtheorem{Co}{Corollary}{\alph{enumi}}
\newenvironment{corollary}{\begin{Co}\hspace{-0.2cm}: }{\end{Co}}

\newtheorem{As}{Assumption}{\alph{enumi}}

\newtheorem{Le}{Lemma}{\alph{enumi}}
\newenvironment{lemma}{\begin{Le}\hspace{-0.2cm}: }{\end{Le}}

\newtheorem{Fo}{Folgerung}{\alph{enumi}}

\newtheorem{De}{Definition}{\alph{enumi}}

\newtheorem{Be}{Bemerkung}{\alph{enumi}}

\newtheorem{Ex}{Example}{\alph{enumi}}
\newenvironment{example}{\begin{Ex}\hspace{-0.2cm}:\rm }{\end{Ex}}

\newcommand{\ol}{\overline}
\newcommand{\bee}{\begin{equation}}
\newcommand{\ee}{\end{equation}}

\newcommand{\bea}{\begin{eqnarray}}
\newcommand{\ea}{\end{eqnarray}}

\newcommand{\ve}{\varepsilon}

\DeclareMathOperator \diver {div}

\begin{document} 
\begin{center}
{\Large{\sc On the Dirichlet problem for the
prescribed mean curvature equation over general domains}} \\[2cm]
{\large{\sc Matthias Bergner}}\\[1cm]
{\small\bf Abstract}\\[0.4cm]
\begin{minipage}[c][3cm][l]{12cm}
{\small 
We study and solve the Dirichlet problem for 
graphs of prescribed mean curvature in $\mathbb R^{n+1}$ 
over general domains $\Omega$ without requiring a mean convexity assumption. 
By using pieces of nodoids as barriers we first
give sufficient conditions for the solvability in case
of zero boundary values. Applying a result by Schulz
and Williams we can then also solve the Dirichlet problem
for boundary values satisfying a Lipschitz condition.} \end{minipage} 
\end{center}\noindent \\[0cm]
\footnote{\noindent Matthias Bergner\\
Universit\"at Ulm, Fakult\"at f\"ur Mathematik und Wirtschaftswissenschaften, Institut f\"ur Analysis\\
Helmholtzstr. 18, D-89069 Ulm, Germany\\
e-mail: matthias.bergner@uni-ulm.de}

\subsection{Introduction}
In this paper we study and solve the Dirichlet problem for
$n$-dimensional graphs of prescribed mean curvature in $\mathbb R^{n+1}$:
Given a domain $\Omega\subset\mathbb R^n$ and 
Dirichlet boundary values $g\in C^0(\partial\Omega,\mathbb R)$
we want to find a solution $f\in C^2(\Omega,\mathbb R)\cap C^0(\ol\Omega,\mathbb R)$ of
\bee
\diver\frac{\nabla f}{\sqrt{1+|\nabla f|^2}}
=n H(x,f) \quad \mbox{in}\; \Omega \label{l1} \; , \;
f=g \quad \mbox{on}\; \partial\Omega \;  .
\ee
The given function $H:\ol\Omega\times\mathbb R\to\mathbb R$
is called the prescribed mean curvature.
At each point $x\in\Omega$ the geometric mean curvature of
the graph $f$, defined as the average of the principal curvatures,
is equal to the value $H(x,f(x))$, thus a solution 
$f$ is also called a graph of prescribed mean curvature $H$. \\ \\
For the minimal surface case, i.e. $H\equiv 0$, it is known that the 
mean convexity of the domain $\Omega$ yields a necessary and sufficient condition
for the Dirichlet problem to be solvable for all Dirichlet boundary values
(see \cite{jenkins}). Here, mean convexity means that $\hat H(x)\geq 0$
for the mean curvature of $\partial\Omega$ w.r.t. the inner normal.
For the prescribed mean curvature case, 
a stronger assumption is needed on the domain $\Omega$
in order to solve the boundary value problem for all Dirichlet boundary values $g$.
A necessary condition on the domain $\Omega$ and the prescribed
mean curvature $H$ is
\bee\label{a1}
|H(x,z)|\leq \frac{n-1}{n} \hat H(x) \quad 
\mbox{for}\; (x,z)\in\partial\Omega\times\mathbb R 
\ee
(see \cite[Corollary 14.13]{gilbarg}).
Additionally requiring a smallness condition
on $H$ implying the existence of a $C^0$-estimate (such as \cite[(10.32)]{gilbarg})
Gilbarg and Trudinger \cite[Theorem 16.9]{gilbarg}
could then solve the Dirichlet problem in case $H=H(x)$.  \\ \\
It is now a natural question to ask if we can 
relax the mean convexity assumption (\ref{a1}) if we only consider
certain boundary values, for example zero boundary values. This is indeed possible,
as our first existence result demonstates.
\begin{theorem}\label{t1}
Assumptions:
\begin{itemize}
\item[a)] Let the bounded $C^{2+\alpha}$-domain $\Omega\subset\mathbb R^n$
satisfy a uniform exterior sphere condition of radius $r>0$
and be included in the annulus $\{x\in\mathbb R^n\; : \; r<|x|<r+d\}$
for some constant $d>0$. 
\item[b)] Let the prescribed mean curvature $H=H(x,z)\in C^{1+\alpha}
(\ol\Omega\times\mathbb R,\mathbb R)$ 
satisfy $H_z\geq 0$ and the smallness assumption
\bee\label{l70}
h:=\sup_{x\in\Omega} |H(x,0)|<\frac{2 (2r)^{n-1}}{(2r+d)^n-(2r)^n} \; . 
\ee
\end{itemize}
Then the Dirichlet problem (\ref{l1}) has a unique
solution $f\in C^{2+\alpha}(\ol\Omega,\mathbb R)$
for zero boundary values.
\end{theorem}
For dimension $n=2$ and constant mean curvature, 
similar existence theorems, again for zero boundary values,
can be found in \cite{lopez2}, \cite{lopez} or \cite{ripoll1}. 
Note that Theorem \ref{t1} can be applied in particular to the
annulus $\Omega:=\{x\in\mathbb R^n\, : \, r<|x|<r+d\}$ which
does not satisfy the mean convexity assumption (\ref{a1}). 
Given any bounded $C^2$-domain $\Omega$ we can find
constants $r>0$ and $d>0$ such that assumption a) of Theorem \ref{t1}
is satisfied for a suitable translation of $\Omega$. \\ \\
The uniqueness part of Theorem \ref{t1} follows directly from
the assumption $H_z\geq 0$ together with the maximum principle.
However, $H_z\geq 0$ is not only needed for the uniqueness but
also for the existence of a solution. More precisely, it
is needed to obtain a global gradient estimate for solutions
of Dirichlet problem (\ref{l1}) (see Theorem \ref{tn}). \\ \\
The smallness condition (\ref{l70}) is required for two
reasons: first to obtain an estimate of the $C^0$-norm
of the solution and secondly to obtain a boundary gradient
estimate (see Theorem \ref{t3}). Other smallness conditions
assuring the existence a $C^0$-estimate are given in \cite{gilbarg}, such as
\bee\label{l240}
h<\Big(\frac{\omega_n}{|\Omega|}\Big)^{1/n} \; . 
\ee
These two assumptions (\ref{l70}) and are (\ref{l240}) quite different as they involve different quantities:
(\ref{l70}) contains the numbers $r$ and $d$ while
(\ref{l240}) contains the volume $|\Omega|$ of $\Omega$.
Additionally, assumption (\ref{l240}) does not imply a boundary gradient estimate
while (\ref{l70}) does. We also want to remark that there are certain
domains for which (\ref{l70}) is satified and not (\ref{l240})
while for certain other domains (\ref{l240}) is satisfied but (\ref{l70}) is not. \\ \\
Note that some kind of smallness assumption on $h$
in Theorem \ref{t1} is needed since there exists
the following necessary condition:
If there exists a graph of constant mean curvature
$h>0$ over a domain $\Omega$ containing a disc
of radius $\varrho>0$, then we have necessarily 
$h\leq \frac{1}{\varrho}$. This follows 
from a comparision with spherical caps of constant mean curvature $\frac{1}{\varrho}$
together with the maximum principle.
Consequently, the smallness condition on $h$ 
in Theorem \ref{t1} cannot solely depend 
on the radius $r$ of the exterior sphere condition. \\
Furthermore, the smallness condition on $h$ also cannot solely depend
on the diameter of the domain: Consider the annulus 
$\Omega=\{x\in\mathbb R^n \; : \; \ve<|x|<1\}$ for some
$0<\ve<1$ with $\mbox{diam}(\Omega)=2$.
In Lemma \ref{lemma_neu} we show that a graph of constant mean
curvature $h>0$ having zero boundary values does not exist
if one chooses $\ve>0$ sufficiently small. \\ \\
Theorem \ref{t1} specifically applies to convex domains.
Note that a convex domain satisfies
a uniform exterior sphere condition of any radius $r>0$.
By letting $r\to+\infty$, we then obtain the following
corollary, which for dimension $n=2$ and constant mean curvature
can also be found in \cite[Corollary 3]{ripoll1} or \cite[Theorem 1.4]{lopez2}.
\begin{corollary}\label{c1}
Let a bounded convex $C^{2+\alpha}$-domain $\Omega\subset\mathbb R^n$
be given such that $\ol\Omega$ is included 
within the strip $\{x\in\mathbb R^n\; | \; 0<x_1<d\}$ of width $d>0$.
Let the prescribed mean curvature $H\in C^{1+\alpha}(\ol\Omega\times\mathbb R,\mathbb R)$
satisfy $H_z\geq 0$ as well as
$$ h:=\sup\limits_{\Omega}|H(x,0)|<\frac{2}{nd} \; . $$
Then the Dirichlet problem (\ref{l1}) has a unique
solution $f\in C^{2+\alpha}(\ol\Omega,\mathbb R)$
for zero boundary values.
\end{corollary}
Note that in Corollary \ref{c1} the diameter
of the domain $\Omega$ can be arbitrarily large, while
in Theorem~\ref{t1} the diameter is bounded by $2(r+d)$. 
Additionally, we can choose the volume $|\Omega|$ 
of the domain $\Omega$ arbitrarily large
so that the smallness assumption (\ref{l240}) will not be satisfied. \\ \\
In case of arbitrary boundary values $g$,
Williams \cite{williams} could show that 
the Dirichlet problem (\ref{l1}) for $H\equiv 0$ is still solvable
over domains not being mean convex domains,
if one requires certain smallness assumptions on $g$. More precisely he showed: 
For any Lipschitz constant $0\leq L<\frac{1}{\sqrt{n-1}}$
there exists some $\ve=\ve(L,\Omega)>0$ such that the Dirichlet
problem (\ref{l1}) is solvable for the minimal surface equation if the
boundary values $g$ satisfy
\bee\label{la}
|g(x)-g(y)|\leq L|x-y|\quad \mbox{for}\; x,y\in\partial\Omega
\quad \mbox{and} \quad |g(x)|\leq \ve \quad \mbox{for}\; x\in\partial\Omega \; . 
\ee
Note that the boundary values are only required to be Lipschitz 
continuous and they are not of class $C^{2+\alpha}$. Hence, also the
solution will be at most Lipschitz continuous up to the boundary.
For the proof Williams first considers weak solutions of the minimal
surface equation. Constructing suitable barriers he then
shows that these weak solutions are continuous up to the
boundary and that the Dirichlet boundary values are attained. \\ \\
Schulz and Williams \cite{schulz}
generalised the result of Williams \cite{williams} 
from the minimal surface case to the prescribed mean curvature case $H=H(x,z)$.
However, two more assumptions are needed there:
As in Theorem \ref{t1}, the prescribed mean curvature function $H$ must
satisfy the monotonocity assumption $H_z\geq 0$.
This assumption is needed for the existence of weak
solutions (see \cite{miranda}).
Moreover, they require the existence of an initial solution 
$f_0\in C^2(\Omega,\mathbb R)\cap C^1(\ol\Omega,\mathbb R)$
for Dirichlet boundary values $g_0$, which must be Lipschitz
continuous with a Lipschitz constant smaller than $\frac{1}{\sqrt{n-1}}$. \\ \\
Using our solution of Theorem \ref{t1}
and Corollary \ref{c1} as an initial solution with zero
boundary values, we can apply the result of Schulz and Williams
to solve the Dirichlet problem for 
Lipschitz continuous boundary values as well:
\begin{theorem}\label{t2}
Let the assumptions of Theorem \ref{t1} or Corollary
\ref{c1} be satisfied. Then for any Lipschitz
constant $0\leq L<\frac{1}{\sqrt{n-1}}$ there
exists some $\ve=\ve(\Omega,H,L)>0$ such that
the Dirichlet problem (\ref{l1}) has a
solution $f\in C^{2+\alpha}(\Omega,\mathbb R)\cap C^{0}(\ol\Omega,\mathbb R)$
for all Lipschitz continuous boundary values 
$g:\partial\Omega\to\mathbb R$ satisfying assumption (\ref{la}).
\end{theorem}
As demonstrated in \cite{schulz}, the smallness assumption
on the Lipschitz constant $L$ is sharp. 
In case of the minimal surface equation, 
Theorem \ref{t2} will be false for any
Lipschitz constant $L>\frac{1}{\sqrt{n-1}}$ and any 
domain $\Omega$ which is not mean convex (see \cite[Theorem 4]{williams}). \\ \\
This paper is organized as follows:
In Section 2 we first we show that solutions
satisfy a height as well as a boundary gradient estimate.
As barriers we use a piece of a rotationally symmetric
surface of constant mean curvature $h$, a so-called Delaunay nodoid. This surface
is constructed in Proposition \ref{p1} by solving
an ordinary differential equation. There we need a smallness assumption
on $h$ corresponding to assumption (\ref{l70}) of Theorem \ref{t1}.
In Section 3 we first give a global gradient estimate
in terms of the boundary gradient (see Theorem \ref{tn}). 
The monotonocity assumption $H_z\geq 0$ plays an important role there.
We then give the proof of Theorem \ref{t1} and Corollary \ref{c1}
using the Leray-Schauder method from \cite{gilbarg}.


\subsection{Estimates of the height and the boundary gradient}
To obtain a priori $C^0$ estimates as well as
boundary gradient estimates for solutions of problem (\ref{l1}),
it is essential to have certain
super and subsolutions at hand serving us upper and lower barriers.
In this paper we will use a rotationally symmetric surface
of constant mean curvature $h$, a so-called Denaunay surface as barrier.
For $h=0$ we have the family of catenoids and for $h\neq 0$
a family consisting of two types of surfaces: 
the embedded unduloids and the immersed nodoids
(see \cite{hsiang}; \cite{kenmotsu} for $n=2$).
We will now construct a piece of the $n$-dimensional 
catenoid (if $h=0$) and $n$-dimensional nodoid (if $h\neq 0$)
which is given as a graph defined over the annulus
$$ \{x\in\mathbb R^n\; | \; r\leq |x|\leq R\} \; . $$
It can be represented almost explicitely 
by solving a second order ordinary differential equation.
\begin{proposition}\label{p1}
Let the numbers $r>0$, $h\geq 0$ and $R>r$ be given satisfying
\bee\label{l6}
h<\frac{{2 (2r)^{n-1}}}{(R+r)^n-(2r)^n} \; . 
\ee
Then there exists a function $p\in C^2([r,R],[0,+\infty))$
with $p(r)=0$ and $p(t)>0$ for $t\in (r,R]$ such that the rotationally symmetric graph
$f(x):=p(|x|)$ defined on the annulus $r\leq |x|\leq R$ has constant mean curvature $-h$.
Furthermore, there exists some $t_0\in (r,R]$ such that
$p(t)$ is increasing for $t\in [r,t_0]$ and decreasing
for $t\in [t_0,R]$.
\end{proposition}  
{\it Proof:}
\begin{itemize}
\item[1.)] Inserting $p(|x|)=f(x)$ into
the mean curvature equation
$$ \diver \frac{\nabla f}{\sqrt{1+|\nabla f|^2}}=-n h $$
we obtain for $p$ the second order differential equation
$$ \frac{p''}{(1+p'^2)^{\frac{3}{2}}}+\frac{(n-1) p'}{t (1+p'^2)^{\frac{1}{2}}}=-n h \; . $$
Multiplying this equation by $t^{n-1}$ and integrating
this yields the first order differential equation
\bee\label{l2}
\frac{t^{n-1}p'}{\sqrt{1+p'^2}}=c-ht^n 
\ee
where $c\in\mathbb R$ is some integration constant serving as a parameter. 
We focus here on the case $c>0$, corresponding to the
choice of a nodoid. The case $c=0$ yields a sphere and $c<0$ 
an unduloid. Solving equation (\ref{l2}) for $p'$ we obtain
\bee\label{l3}
p'(t)=\frac{c-ht^n}{\sqrt{t^{2n-2}-(c-h t^n)^2}} \; . 
\ee
Clearly, (\ref{l3}) is only well defined for those $t\in (0,+\infty)$
for which the term under the root in the denominator is positive.
We will later determine for which $t$ this is the case.
Integrating (\ref{l3}) we can now define
\bee\label{l250}
p(t):=\int\limits_r^t \frac{c-hs^n}{\sqrt{s^{2n-2}-(c-h s^n)^2}} ds  
\ee
with $p(r)=0$.
\item[2.)] Let us first study the case $h=0$. 
The denominator of (\ref{l3}) has exactly one zero
$a>0$ given as solution of $a^{n-1}=c$ and $p'(t)$ is defined
for all $t\in (a,+\infty)$. For the integral
(\ref{l250}) to be defined, we need to have that $r\in (a,+\infty)$, which is equivalent
to $c<r^{n-1}$. For example, we can set $c:=\frac{1}{2} r^{n-1}$.
The function $p(t)$ is now defined for
all $t\in [r,+\infty)$ and also $p'(t)>0$ for all $t\in [r,+\infty)$.
The claim of the proposition now follows with $t_0=R$.
\item[3.)] In case $h>0$, the denominator of (\ref{l3}) 
has precisely two positive zeros
$0<a<b$ given as solutions of the equations
$$ h a^n+a^{n-1}=c \quad , \quad hb^n-b^{n-1}=c \; . $$
Now $p'(t)$ is defined for all $t\in (a,b)$
and formally we have $p'(a)=+\infty$, $p'(b)=-\infty$.
Note that for
$$ t_0:=\Big(c\, h^{-1}\Big)^{\frac{1}{n}}\in (a,b) $$
we have 
$$ p'(t_0)=0 \quad , \quad  p'(t)>0 \quad \mbox{for}\; t\in (a,t_0)
\quad \mbox{and} \quad p'(t)<0 \quad \mbox{for}\; t\in (t_0,b) \; , $$
as desired. Now for the integral (\ref{l250})
to be defined, we need to have $a<r<t_0$, which is
equivalent to restricting the parameter $c$ such that
\bee\label{l270}
h r^n<c<h r^n+r^{n-1} \; . 
\ee
We then obtain $p\in C^2([r,b),\mathbb R)$.
\item[4.)] We will now show the inequality 
\bee\label{l5}
p'(t_0-s)>|p'(t_0+s)| \quad \mbox{for all}\; s\in (0,t_0-a) \; .
\ee
Together with $p(r)=0$ this will yield $p(t)>0$ for all $t\in (r,r+2(t_0-r)]$.
In fact, after some computation (\ref{l5})
turns out to be equivalent to
$$ q(t_0-s)+q(t_0+s)>0 \quad \mbox{for}\; s\in (0,t_0-a) $$
for the function $q(t):=(c-ht^n) t^{1-n}=c t^{1-n}-h t$.
This however is a direct consequence of the inequality
$$ c(t_0+s)^{1-n}+c(t_0-s)^{1-n}>2h t_0 $$
which holds for all $s\in (0,t_0)$, proving (\ref{l5}).
\item[5.)] We now set
$$ R'=R'(c):=r+2(t_0-r)=2t_0-r=2\Big(c h^{-1}\Big)^{\frac{1}{n}}-r<b \; . $$
From 4.) we conclude the positivity $p(t)>0$ for all $t\in (r,R']$. 
Keeping in mind the restriction (\ref{l270})
on $c$ we obtain the limit
$$ R'(c)\to 2\big(r^n+ h^{-1} r^{n-1}\Big)^{\frac{1}{n}}-r
=2r\Big(1+h^{-1} r^{-1}\Big)^{\frac{1}{n}}-r $$
if we let $c\to h r^n+r^{n-1}$.
This proves the claim of the proposition whenever
$$ R<2r\Big(1+h^{-1} r^{-1}\Big)^{\frac{1}{n}}-r $$
is satisfied. An easy computation, however, asserts that
this inequality is indeed equivalent to assumption (\ref{l6}) . \hfill $\Box$ 
\end{itemize}
The following picture shows the graph of the function $p(t)$
for $n=2$, $h=\frac{1}{3}$, $a=1$ and $b=4$.
\psfrag{r}{$r$}
\psfrag{a}{$a$}
\psfrag{R}{$R$}
\psfrag{b}{$b$}
\psfrag{t0}{$t_0$}
\includegraphics[scale=0.88]{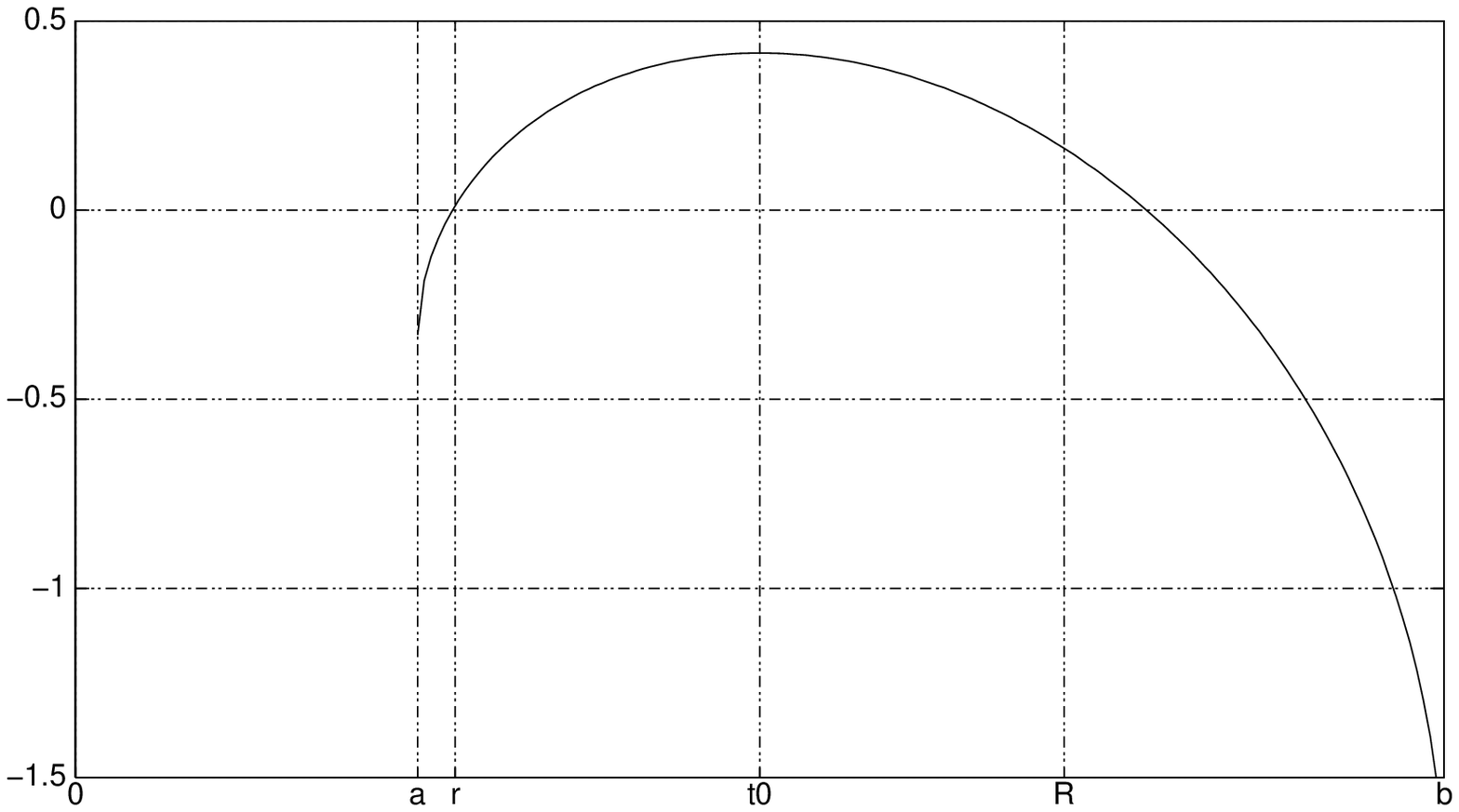}\\ \noindent 
{\it Remarks:} 
\begin{itemize}
\item[a)] For $h=0$ and $n=2$ the function
$p(t)$ has the explicit form 
$p(t)=c\, \mbox{arcosh}(t/c)$, i.e. the well known catenary.
If either $h>0$ or $n\geq 3$ the function $p(t)$ 
can only be represented by the elliptic integral given in the 
proof of Proposition \ref{p1}.
\item[b)] In the case $h=0$ we obtain the $n$-dimensional
catenoid, a rotationally symmetric minimal surface.
The generating function is defined for all $t\in [r,+\infty)$.
In case $n=2$ we have  $p(t)\to\infty$ as $t\to\infty$. 
However, for $n\geq 3$ the function
$p(t)$ is uniformly bounded by some constant.
\item[c)] In case $h>0$, the maximal domain of definition 
of the function $p(t)$ is the interval $(a,b)$. In case $n=2$ one can show that the
length $b-a$ of this interval is given by $b-a=\frac{1}{h}$,
in particular the length does not depend on the parameter $c$.
This is no longer the case for dimension $n\geq 3$
where $b-a$ depends on both $h$ and $c$.
\end{itemize}
At this point let us prove the following nonexistence result 
which we already claimed in the introduction.
\begin{lemma}\label{lemma_neu}
For $0<\ve<1$ consider the annulus $\Omega:=\{x\in\mathbb R^n\, : \, \ve<|x|<1\}$.
Then given any constant $h>0$ there exists some $\ve=\ve(h)\in (0,1)$, such that
a graph $f\in C^2(\Omega,\mathbb R)\cap C^0(\overline\Omega,\mathbb R)$
of constant mean $h$ with zero boundary values does not exist. 
\end{lemma}
{\it Proof:}
We will show that such a graph of constant mean curvature $-h$ does not exist
for sufficiently small $\ve>0$. By a reflection argument, then a graph of
constant mean curvature $h$ does not exist either.
Assume to the contrary that a graph $f=f_\ve$ does exist for each $\ve>0$.
Because $f_\ve$ has constant mean curvature $-h<0$ and zero boundary values, 
the maximum principle yields $f_\ve(x)\geq 0$ for $x\in\Omega$.
Now note that the domain $\Omega$ and the boundary values of $f_\ve$ are rotationally symmetric. 
Hence, the solution $f_\ve$ is also rotationally symmetric, following from the uniqueness of
the Dirichlet problem. But then we can write $f_\ve(x)=p_\ve(|x|)$ where $p_\ve(t)$ satisfies
$p_\ve(t)\geq 0$ for $t\in [\ve,1]$ and $p_\ve(\ve)=p_\ve(1)=0$. From (\ref{l3}) we conclude
$$ p_\ve(t)=\int\limits_\ve^t \frac{c-hs^n}{\sqrt{s^{2n-2}-(c-h s^n)^2}} ds $$
where $c=c(\ve)\in\mathbb R$ is a suitable constant. 
We set $k:=c-h\ve^n$ and claim $k\geq 0$. Otherwise $p_\ve'(t)<0$ would hold for 
all $t\in (\ve,1)$, contradicting $p_\ve(\ve)=p_\ve(1)=0$.
Note that the expression under
the root must be nonnegative for all $s\in [\ve,1]$, in particular for $s=\ve$ we get
$$ \ve^{2n-2}-(c-h \ve^n)^2=\ve^{2n-2}-k^2\geq 0  $$
or equivalently $k^{-2} \ve^{2n-2}\geq 1$. For any $t\in [\ve,1]$ we now estimate
\bea
p_\ve(t)&=&\int_\ve^t \frac{c-hs^n}{\sqrt{s^{2n-2}-(c-h s^n)^2}} ds 
\leq\int_\ve^t \frac{c-h\ve^n}{\sqrt{s^{2n-2}-(c-h \ve^n)^2}} ds  \nonumber \\
&=&\int_1^{t/\ve} \frac{k}{\sqrt{(\ve \tau)^{2n-2}-k^2}} \ve d\tau
=\ve \int_1^{t/\ve} \frac{1}{\sqrt{k^{-2}\ve^{2n-2} \tau^{2n-2}-1}}d\tau  \nonumber \\
&\leq&\ve \int_1^{1/\ve}\frac{1}{\sqrt{\tau^{2n-2}-1}}d\tau \nonumber \; . 
\ea
In case of dimension $n\geq 3$ we conclude that $\lim\limits_{\ve\to 0} p_\ve(t)=0$,
which follows from
$$ \int_1^{+\infty}\frac{1}{\sqrt{\tau^{2n-2}-1}}d\tau<+\infty  $$
for dimension $n\geq 3$. In case of $n=2$, the above integral is infinite. However, the explicit computation
$$ p_\ve(t)\leq\ve \int_1^{1/\ve}\frac{1}{\sqrt{\tau^2-1}} d\tau
=\ve \Big[\mbox{arcosh}(t)\Big]_1^{1/\ve}=\ve\mbox{arcosh}(1/\ve) $$
again shows $\lim\limits_{\ve\to 0} p_\ve(t)=0$.
This implies that the family $f_\ve(x)=p_\ve(|x|)$ converges uniformly 
to $f_0(x)\equiv 0$ on every compact subset of $\{x\in\mathbb R^n \, : \, 0<|x|\leq 1\}$.
Then, after extracting some subsequence, all first and second derivatives of $f_\ve$ will converge to zero
by interior gradient estimates for the constant mean curvature equation.
Hence, also the mean curvature of $f_\ve$ must converge to zero.
This yields a contradiction as the mean curvature of $f_\ve$ is $-h$ for
each $\ve>0$. \hfill $\Box$\\ \\
We can now show the a priori estimates
of the height and boundary gradient of solutions of (\ref{l1}).
\begin{theorem}\label{t3}
Assumptions:
\begin{itemize}
\item[a)] Let the bounded $C^2$-domain $\Omega\subset\mathbb R^n$
satisfy a uniform exterior sphere condition of radius $r>0$
and be included in the annulus $\{x\in\mathbb R^n\; : \; r<|x|<r+d\}$
for some constant $d>0$. 
\item[b)] Let the prescribed mean curvature $H=H(x,z)\in C^1
(\ol\Omega\times\mathbb R,\mathbb R)$ 
satisfy $H_z\geq 0$ in $\Omega\times\mathbb R$ 
as well as the smallness assumption $|H(x,0)|\leq h$
for some constant
$$ h<\frac{2 (2r)^{n-1}}{(2r+d)^n-(2r)^n} \; . $$
\item[c)] Let $f\in C^2(\ol\Omega,\mathbb R)$
be a solution of problem (\ref{l1}) for 
zero boundary values.
\end{itemize}
Then there exists a constant $C=C(h,r,d)$ such that $f$
satisfies the estimates
$$ ||f||_{C^0(\Omega)}\leq C \quad \mbox{and} \quad 
\sup\limits_{\partial\Omega} |\nabla f(x)|\leq C \; . $$
\end{theorem}
{\it Proof:}
\begin{itemize}
\item[1.)] We first show the $C^0$-estimate.
Because of $\Omega\subset \{x\in\mathbb R^n \, : \, r<|x|<r+d\}$
the rotationally symmetric graph $\eta(x):=p(|x|)$
is well defined and has constant mean curvature $-h$.
Here, $p(t)$ is the function defined by Proposition \ref{p1} for $R:=r+d$.
From $|H(x,0)|\leq h$ together with $H_z\geq 0$ we conclude
\bee\label{l120}
H(x,z)\geq -h \quad \mbox{for} \; x\in\Omega \; , \, z\geq 0
\quad \mbox{and} \quad 
H(x,z)\leq h \quad \mbox{for} \; x\in\Omega \; , \, z\leq 0 \; .
\ee
We now choose $c\geq 0$ minimal such that $f(x)\leq \eta(x)+c$
holds in $\ol\Omega$. We claim that $c=0$. Otherwise
there would be a point $x_0\in\Omega$ with $f(x_0)=\eta(x_0)+c>0$.
From (\ref{l120}) together with the strong maximum principle
we then would have $f(x)\equiv \eta(x)+c$ in $\Omega$, 
contradicting $f(x)=0$ on $\partial\Omega$.
Hence we have shown $f(x)\leq \eta(x)$ in $\Omega$. Similary, we obtain
$f(x)\geq -\eta(x)$. Combining these estimates we have
$$ ||f||_{C^0(\Omega)}=\sup\limits_{\Omega} |f(x)|\leq
\sup\limits_{\Omega} |\eta(x)|
\leq \sup\limits_{r\leq t\leq r+d} |p(t)|=p(t_0)=: C_1 \; . $$
Here, $t_0$ defined by Proposition \ref{p1} is
the argument for which the function $p$ achieves its maximum.
Note that $p$ only depends on $r,d$ and $h$ and hence $C_1=C_1(r,d,h)$.
\item[2.)] Given some point $x_0\in\partial\Omega$ we show the 
boundary gradient estimate at $x_0$. 
Since $\Omega$ satisfies a uniform exterior sphere condition of radius $r$,
we may assume that
$$ \Omega\cap B_r(0)=\emptyset \quad \mbox{and}
\quad x_0\in\partial B_r(0)\cap\partial\Omega $$
holds after a suitable translation. 
We define the annulus $U:=\{x\in\mathbb R^n\, :\, r<|x|<t_0\}$
and consider the graph
$$ \eta\in C^2(\ol U,\mathbb R) \quad , \quad \eta(x):=p(|x|) \quad \mbox{for}\;
x\in\ol U\; . $$
From $f(x)=0$ on $\partial\Omega$ together with
$f(x)\leq p(t_0)=\eta(x)$ for $|x|=t_0$
we conclude $f(x)\leq \eta(x)$ on $\partial(\Omega\cap U)$.
As in part 1.), the maximum principle gives
$f(x)\leq \eta(x)$ in $\Omega\cap U$ as well as
$f(x)\geq -\eta(x)$ in $\Omega\cap U$.
From $x_0\in\partial(\Omega\cap U)$
and $f(x_0)=\eta(x_0)$ we obtain
$$ |\nabla f(x_0)|=\Big|\frac{\partial}{\partial \nu} f(x_0)\Big|
\leq \Big|\frac{\partial}{\partial \nu} \eta(x_0)\Big|=|p'(r)|=: C_2 \; , $$
where $\nu$ is the outward normal to $\partial\Omega$ at $x_0$. \hfill $\Box$
\end{itemize}
{\it Remark:} A closer inspection of the proof shows that
Theorem \ref{t3} also holds without the assumption
$H_z\geq 0$ if one requires $|H(x,z)|\leq h$ 
in $\Omega\times\mathbb R$ instead of $|H(x,0)|\leq h$ in $\Omega$.
However, we will essentially need the assumption $H_z\geq 0$ in the next section to
prove a global gradient estimate.


\subsection{Global gradient estimate and the proof of Theorem \ref{t1}}
In the previous section we have shown 
a $C^0$-estimate together with a boundary
gradient estimate, thus we can assume
\bee\label{l55}
|f(x)|\leq M \quad \mbox{in}\; \ol\Omega 
\ee
for a given solution $f\in C^{2+\alpha}(\ol\Omega,\mathbb R)$
of problem (\ref{l1}).
It now remains to establish a global gradient estimate
in terms of the $C^0$-norm and the boundary gradient.
Such a global gradient estimate is derived in \cite[Theorem 15.2]{gilbarg}
for a fairly large class of quasilinear elliptic equations.
This includes the prescribed mean curvature equation in case of $H=H(x)$,
as verified in example (ii) after \cite[Theorem 15.2]{gilbarg}. We will show that 
\cite[Theorem 15.2]{gilbarg} continues to hold in case $H=H(x,z)$, if we assume the
monotonocity condition $H_z\geq 0$. 
Let us first write the prescribed mean curvature equation in the form
$$ \triangle f-\sum_{i,j=1}^n \frac{\partial_i f\partial_j f}{1+|\nabla f|^2} \partial_{ij} f- nH(x,f) \sqrt{1+|\nabla f|^2}=0 \; . $$
Now quantities $\alpha,\beta,\gamma$ are defined by \cite[(15.27)]{gilbarg}, which in our case are
\begin{eqnarray}
&&\alpha=-1+\frac{1}{1+|p|^2} \quad , \quad \beta=\frac{n H(x,z) \sqrt{1+|p|^2}}{|p|^2}  \nonumber \\
&&\gamma=-n \frac{(1+|p|^2)^{3/2}}{|p|^2}
\Big[H_z(x,z)+\sum_{i=1}^n \frac{p_i}{|p|^2}  H_{x_i}(x,z)\Big] \quad 
\mbox{for}\; x\in\Omega \; , \; |z|\leq M \; , \; p\in\mathbb R^n \nonumber 
\end{eqnarray}
(compare with example (ii) in chapter 15.2 of \cite{gilbarg}).
We now compute the limits
\begin{eqnarray}
&&a:=\limsup\limits_{|p|\to\infty} \alpha=-1 \quad , \quad 
b:=\limsup\limits_{|p|\to\infty} \beta\leq n \sup\limits_{\Omega\times [-M,M]} |H(x,z)| \nonumber \\
&& c:=\limsup\limits_{|p|\to\infty} \gamma\leq n \sup\limits_{\Omega\times [-M,M]} |\nabla H(x,z)| \label{l27}
\end{eqnarray}
using $H_z\geq 0$ for the last limit. Because of $a=-1$ together with $b,c<+\infty$ we may apply 
\cite[Theorem 15.2]{gilbarg} to obtain
\begin{theorem}\label{tn}
Let the prescribed mean curvature 
$H\in C^1(\overline\Omega\times\mathbb R,\mathbb R)$ satisfy
$$ H_z(x,z)\geq 0 \quad , \quad |H(x,z)|+|\nabla H(x,z)|\leq h_0 \quad \mbox{for}\; x\in\Omega \; , \; |z|\leq M \; . $$
Let $f\in C^2(\overline\Omega,\mathbb R)$ be a solution Dirichlet problem (\ref{l1})
satisfying $||f||_{C^0(\Omega)}\leq M$.
Then the estimate
$$ \sup\limits_{x\in\Omega} |\nabla f(x)|\leq C $$
holds with a constant $C$ depending only on $n$, $h_0$, $M$, $\Omega$
and $\sup_{\partial\Omega} |\nabla f|$.
\end{theorem}
{\it Remark:} If we do not assume $H_z\geq 0$, then we will obtain
$c=+\infty$ in (\ref{l27}) and \cite[Theorem 15.2]{gilbarg} will not be applicable. 
In fact, the following example shows that a gradient estimate is false if one
does not require $H_z\geq 0$.
\begin{example}
Given some parameter $\ve>0$ let $\beta(z):=z^3+\ve z$ for $z\in I:=[-1,1]$.
Noting $\beta'(z)=3z^2+\ve>0$ in $I$, there exists a smooth inverse 
$\beta^{-1}:I\to\mathbb R$. From $\beta(-1)\leq -1$ and $\beta(1)\geq 1$ we 
conclude $\beta^{-1}:I\to I$. We now consider the one-dimensional graph
$f_\ve(x):=\beta^{-1}(x)$ for $x\in I$ with its parametrisation $X(x)=(x,f_\ve(x))$.
Substituting $z=f_\ve(x)$ we obtain the reparametrisation $\tilde X(z)=(\beta(z),z)$
and we can compute the curvature $H=H(z)$ by
$$ H(z):=H_\ve(z)=-\frac{\beta''}{\Big(1+(\beta')^2\Big)^{3/2}}=-\frac{6z}{\Big(1+(3z^2+\ve)^2\Big)^{3/2}} \; . $$
Hence, $f_\ve$ is a graph of prescribed mean curvature $H_\ve(z)$.
We can find a constant $C$ such that
$$ |H_\ve(z)|+|\nabla H_\ve(z)|\leq C \quad \mbox{for all} \; z\in [-1,1] \; , \; 0< \ve\leq 1 \; . $$
Additionally, we have the $C^0$-estimate and boundary gradient estimate
$$ |f_\ve(x)|\leq 1 \quad \mbox{for}\; x\in I \quad \mbox{and}\quad |\nabla f_\ve(x)|\leq 1 \quad
\mbox{for}\; x\in\partial I=\{-1,1\} \; . $$
However, there is no uniform gradient bound for $f_\ve$ in $I$ because
$$ |\nabla f_\ve(0)|=|f_\ve'(0)|=\frac{1}{|\beta'(0)|}=\frac{1}{\ve}\to \infty \quad \mbox{if}\; \ve\to 0 \; . $$
In this example, all of the assumptions of Theorem \ref{tn} are satisfied except for $H_z\geq 0$.
Even though this example was purely one-dimensional, a generalisation to
higher dimensions $n\geq 2$ is easily possible.
\end{example}
We can now give the\\ \\
{\it Proof of Theorem \ref{c1}:} \\
For $t\in [0,1]$ consider the family of Dirichlet problems
\bee
f\in C^{2+\alpha}(\ol\Omega,\mathbb R) \quad , \quad
\diver\frac{\nabla f}{\sqrt{1+|\nabla f|^2}}
=t\,n H(x,f) \quad \mbox{in}\; \Omega \quad \mbox{and}\quad 
f=0 \quad \mbox{on}\; \partial\Omega\label{l11} \; .
\ee
Let $f$ be such a solution
for some $t\in [0,1]$. By Theorems \ref{t3} and \ref{tn} we have the estimate
$$ ||f||_{C^1(\Omega)}\leq C $$
with some constant $C$ independet of $t$. 
The Leray-Schauder method \cite[Theorem 13.8]{gilbarg}
yields a solution of the Dirichlet problem (\ref{l11}) for each $t\in [0,1]$.
For $t=1$ we obtain the desired solution of (\ref{l1}).  \hfill $\Box$ \\ \\
{\it Proof of Corollary \ref{c1}:} \\
Corollary \ref{c1} is obtained as the limit case
of Theorem \ref{t1} by increasing the radius $r$ of
the exterior sphere condition to infinity.
First, since $\ol\Omega$ is bounded and included within the strip
$\{x\in\mathbb R^n\; : \; 0<x_1<d\}$, after a suitable
translation it will also be included within the
annulus $\{x\in\mathbb R^n\; : \; r<|x|<r+d\}$ for sufficiently
large $r>0$. To show which smallness condition on $h$ is required
in order to apply Theorem \ref{c1} we have to compute the limit
\bee\label{l57}
\lim\limits_{r\to\infty} \frac{2 (2r)^{n-1}}{(2r+d)^n-(2r)^n} \; . 
\ee
To do this, we calculate
$$
\lim\limits_{r\to\infty} \frac{(2r+d)^n-(2r)^n}{2 (2r)^{n-1}} 
=\lim\limits_{r\to\infty} 
\frac{(2r)^n+n (2r)^{n-1} d +O(r^{n-2})-(2r)^n}{2 (2r)^{n-1}}
=\frac{nd}{2} \; . $$
We see that the limit in (\ref{l57}) is equal
to $\frac{2}{n d}$ and hence the smallness condition
$h<\frac{2}{n d}$ is required. Alternatively we could prove 
Corollary \ref{c1} also directly, by proving an analogue
result to Theorem \ref{t3} for convex domains. 
Instead of using the nodoid we would then use a cylinder
as barrier whose axis is lying in the $x_1,\dots,x_n$ hyperplane.
Note that the cylinder $\{x\in\mathbb R^{n+1}\; : x_1^2+\dots+x_n^2=(\frac{d}{2})^2\}$
has constant mean curvature $h=\frac{2}{n d}$, corresponding to the
smallness condition from above. \hfill $\Box$ \\ \\
{\it Remarks:}
\begin{itemize}
\item[a)] Using the methods from \cite{bergner1}, it is also
possible to generalise Theorem \ref{t1} and Corollary \ref{c1}
to the case of prescribed anisotropic mean curvature
$$ \mbox{div}\frac{\nabla f}{\sqrt{1+|\nabla f|^2}}=n H(x,f,N) \quad \mbox{in}\; \Omega \; . $$
Here, the prescribed mean curvature does not only depend on the point $(x,f(x))$ in
space but also on the normal $N(x)$ of the graph.
Within this situation, $H_z\geq 0$ 
can be relaxed to weaker assumption allowing nonuniqueness of solutions.
\item[b)] The results can also be generalised in another direction: 
Define the boundary part
$$ \Gamma_+:=\Big\{x\in\partial\Omega \; : \; 
|H(x,z)|\leq \frac{n-1}{n} \hat H(x) \; \mbox{for all}\; z\in\mathbb R\Big\} $$
where $\hat H(x)$ is the mean curvature of $\partial\Omega$ at $x$ w.r.t.
the inner normal. Now choose a subset $\Gamma\subset \Gamma_+$ such that 
$\mbox{dist}(\Gamma,\partial\Omega\backslash\Gamma_+)>0$.
On $\Gamma$ we can use the standard boundary gradient estimate 
(see \cite[Corollary 14.8]{gilbarg}) and prescribe $C^{2+\alpha}$ boundary values $g$
there. Our boundary gradient estimate of Theorem \ref{t3}, requiring 
zero boundary values, is then only needed on $\partial\Omega\backslash\Gamma$.
This way, Theorem \ref{t1} and Corollary \ref{c1} also hold
for Dirichlet boundary values $g\in C^{2+\alpha}(\partial\Omega,\mathbb R)$
with $g(x)=0$ on $\partial\Omega\backslash\Gamma$ and $|g(x)|\leq \ve$,
where $\ve=\ve(\Omega,\Gamma,H)>0$ is some constant determined
by the height of the nodoid constructed in Proposition \ref{p1}.
\end{itemize}

\end{document}